\documentclass[11pt]{article}
\usepackage[auth-sc]{authblk}
\usepackage{amssymb}
\usepackage{amsthm}
\usepackage{graphicx}
\usepackage{latexsym}
\usepackage{amsmath,amscd}
\usepackage{bm}
\def\part#1{\frac{\partial\phantom{q}}{\partial#1}}

\newenvironment{ex}{\begin{trivlist}\item[]{\bf Example:} }
{\end{trivlist}}

\newenvironment{prf}{\begin{trivlist}\item[]{\bf Proof:} }
{\hfill $\Box$ \end{trivlist}}
\newenvironment{lemprf}{\begin{trivlist}\item[]{\bf Proof:} }
 {\end{trivlist}}
\newtheorem{thm}{Theorem}
\newtheorem{definition}{Definition}
\newtheorem{prp}[thm]{Proposition}
\newtheorem{lemma}[thm]{Lemma}

\newcommand{\lie}[1]{\mathfrak{#1}}
\def\End{\mathop{\rm End}\nolimits}

\def\tr{\mathop{\rm tr}\nolimits}

\def\ad{\mathop{\rm ad}\nolimits}

\def\Sym{\mathop{\rm Sym}\nolimits}

\newcommand{\R}{\mathbf{R}}

\newcommand{\Z}{\mathbf{Z}}

\title{ A note on vanishing theorems}
 \author{Nigel Hitchin}

 \affil  {Mathematical Institute,
Radcliffe Observatory Quarter,
Woodstock Road,
Oxford, OX2 6GG}

\begin{document}
 \maketitle
 \thispagestyle{empty}

\let\oldthefootnote\thefootnote
\renewcommand{\thefootnote}{\fnsymbol{footnote}}
\footnotetext{hitchin@maths.ox.ac.uk}
\let\thefootnote\oldthefootnote
\centerline{\it To the memory of Shoshichi Kobayashi}

\section{Preface}

My own name and that of Professor Kobayashi are now inextricably linked in what has become known as the {\it Hitchin-Kobayashi correspondence}. Despite quizzing the obvious suspects, I do not know who first coined the phrase, but it may be useful to record how, for my part, it was   initially formulated. This was  in a pamphlet ``Non-linear problems in geometry" published by Takeshi Kotake of Tohoku University as a result of a Taniguchi Foundation conference held at Katata, Japan from September 3--8 1979. Each of the participants was asked to produce a list of open problems. I chose three, and the first was the following:
\begin{quotation}

{\bf Problem C-1.} Given a complex manifold $X$ with $c_1(X)=0$ and a positive cohomology class $h\in H^{1,1}(X,\R)$, the Calabi-Yau theorem defines a K\"ahler form $\omega$ which represents $h$ and whose Ricci tensor vanishes. There should be an analogue of this theorem for more general vector bundles than the tangent bundle. More precisely, let $X$ be an algebraic variety with a Hodge metric whose K\"ahler form is $\omega$. Does any holomorphic vector bundle $E$ on $X$ with $c_1(E)=0$ which is {\it stable} relative to the polarization $[\omega]\in H^2(X,\Z)$ have a hermitian metric such that the curvature of the associated connection is orthogonal to $\omega$? 

For curves this is true: orthogonality to $\omega$ means that the curvature vanishes and it is known that any stable bundle has a flat connection. In higher dimensions, the instanton bundles on $P_3$ give examples of such metrics on stable bundles. The K\"ahler metrics of the Calabi-Yau theorem are examples on the tangent bundle. 

Inequalities for characteristic classes and vanishing theorems for holomorphic tensor fields would follow from a positive answer to the question. 
\end{quotation}

For me, the evidence was strong but I had no idea how to prove it and it was only when I suggested the problem the following year to my student Simon Donaldson, and in particular when he gave a moment map interpretation of the equations, that I was convinced it must be true. And of course Donaldson later gave a proof for algebraic surfaces which provided the input   for many examples of his ground-breaking work on gauge theory and four-dimensional topology. 

The vanishing theorems mentioned in the problem were at the time familiar to me from the theorem of Kobayashi and Wu 
\cite{KW} which liberated a basic idea of Yano and Bochner from Riemannian geometry to the realm of gauge theory on a Hermitian manifold and had been important tools both in my thesis on the Dirac operator and in the work on instantons. They play a central role in the later book of Kobayashi \cite{K}.

What follows consists of some observations on vanishing theorems in Riemannian geometry rather than Hermitian geometry. They  date from the mid 1970s but  were never written up, although this point of view was adopted by Arthur Besse in Chapter I of \cite{Bes}. I hope that after 40 years they may still be of some interest.

\section{Introduction}
Laplacian-like operators frequently appear in Riemannian geometry, defined on functions, forms, tensors or spinors. Our aim here is to put all these operators into a general context by defining a Laplacian on any vector bundle $E$ associated to the principal frame bundle  by a representation $\rho$ of  $SO(n)$ or $Spin(n)$. We then define a Laplacian $\Delta$ by means of the so-called Weitzenb\"ock decomposition
$$\Delta=\nabla^*\nabla+tK$$
where $\nabla$ is the covariant derivative of the connection induced on the vector bundle by the Levi-Civita connection, $K$ is a curvature term  and $t$ is a real number. The term $K$ is defined as follows: the Riemann curvature tensor at each point lies in $\Sym^2({\lie so}(n))$. Applying the representation $\rho:{\lie so}(n)\rightarrow \End E$ we get a term in $\Sym^2(\End E)$ and then composition of endomorphisms gives a self-adjoint endomorphism  of $E$. We allow the real number $t$ to vary and as we shall see, different geometrical problems give rise to different values of $t$.

We show that many of the natural Laplacians in geometry such as  the Hodge Laplacian, the Lichnerowicz Laplacian, and the spinor Laplacian are of this form. When the holonomy reduces to a subgroup $H$, it follows directly that  the corresponding spaces of solutions to the Laplace equation decompose into irreducible representations of $H$, a much-used fact in studying manifolds of special holonomy. This applies in particular to the Hodge Laplacian and decomposes the cohomology into irreducible components, generalizing the Hodge decomposition for a K\"ahler manifold. 

Vanishing theorems come from applying Stokes' theorem to the Weitzenb\"ock formula in the case that the curvature term is a positive self-adjoint operator at each point. Here we observe that positivity for {\it all} non-trivial  irreducible representations of $SO(n)$ is equivalent to the manifold having positive curvature operator, which we now know implies that it is  diffeomorphic to a spherical space form \cite{BW}. 

\section{Laplacians}
Let $M$ be an oriented $n$-dimensional Riemannian manifold. We shall write $P$ for the principal bundle of $SO(n)$ or (if $w_2(M)=0$) $Spin(n)$ frames. Given a (unitary) representation $\rho$ of any one of these groups on a vector space $E$ we take the associated vector bundle $P\times_{\rho}E$ which for convenience we shall still denote by the letter $E$. The Levi-Civita connection induces a covariant derivative $\nabla$ on sections of $E$, preserving the inner product. 

The Riemann curvature tensor $R$ is a section of  $\Lambda^2T^*\otimes {\lie so(n)}$ where ${\lie so(n)}$ denotes the bundle associated to the adjoint representation. But using the metric this is isomorphic to $\Lambda^2T^*$, so $R$ is a section of ${\lie so(n)}\otimes {\lie so(n)}$. Using a local orthonormal basis $\{x_a\}$ for ${\lie so(n)}$ we have
$$R=\sum_{a,b} R_{ab}x_a\otimes x_b$$
and the basic symmetries of $ R$ tell us that  $R_{ab}$ is symmetric. 

The representation defines a Lie algebra homomorphism $\rho: {\lie so(n)}\rightarrow \End E$ and then $(\rho\otimes \rho)(R)\in \End E\otimes \End E$. Composing the endomorphisms gives $K$, a section of $\End E$. So
\begin{equation}
K=\sum  R_{ab}\rho(x_a)\rho(x_b).
\label{K}
\end{equation}
Since $\rho(x_a)$ is skew-adjoint and $R_{ab}=R_{ba}$, $K$ is self-adjoint. 

\begin{ex} For the sphere, $R_{ab}$ is the identity matrix and then $K$ is just the quadratic Casimir of the representation $\rho$. More generally, since the square of a skew-adjoint matrix is negative, when $R_{ab}$ is positive, $K$ is negative. 
\end{ex}

\begin{definition} \label{Lap} A Laplacian on $E$ is a differential operator of the form 
$$\Delta=\nabla^*\nabla+tK$$
where $t\in \R$ and $\nabla^*:\Omega^1(E)\rightarrow \Omega^0(E)$ is the formal adjoint of $\nabla: \Omega^0(E)\rightarrow \Omega^1(E).$
\end{definition} 
Clearly $\Delta$ is a second order self-adjoint differential operator. Using Stokes' theorem $\nabla^*\nabla$ can be written as $-\tr\nabla^2$ where $\tr$ denotes the contraction $\tr: E\otimes \Sym^2 T^*\rightarrow E$  using the metric $g$. 

The basic spin representation $\rho$ of $Spin(n)$ on the spinor space $V$ generates all representations in the sense that any irreducible is a subspace of some tensor product $V\otimes \cdots\otimes V$, and we begin by studying the  spinor Laplacian. 
\subsection{The spinor Laplacian} 
The spinor Laplacian is the square of the Dirac operator $D$ acting on sections of the spinor bundle $V$. This is the composition $C\nabla$ where $C:V\otimes T^*\rightarrow V$ is Clifford multiplication $C(\psi\otimes a)=a\cdot\psi$. Then 
$D^2=C^2\nabla^2$ since $C$ is covariant constant. Clifford multiplication satisfies the relation
$$a\cdot b+b \cdot a=-2g(a,b)1,$$
so the symmetric part of $C^2\nabla^2$  gives $-\tr \nabla^2=\nabla^*\nabla$. The skew-symmetric part is defined by the curvature.  Using an orthonormal basis $\{e_i\}$ of $T^*$ the curvature is  a 2-form $\sum K_{ij}e_i\wedge e_j$ with values in the Lie algebra ${\lie so}(n)$ and the action is 
$$\psi\mapsto \sum e_i\cdot e_j\cdot \rho(K_{ij})(\psi).$$ 
 But under the spin representation a skew-symmetric matrix $a_{ij}\in \lie{so}(n)$ maps to the element $-\sum a_{ij}e_i\cdot e_j/4$ in the Clifford algebra, hence the curvature term is in our notation $-4K$ and 
 $$D^2=\nabla^*\nabla-4K.$$
Using the symmetries of the Riemann curvature tensor we have $-4K=s/4$, the Lichnerowicz formula \cite{L1}, nowadays also ascribed to Schr\"odinger. 

An amusing consequence of this approach is to consider the curvature term  for a compact simple group $G$ with metric given by $g=-B$ where $B$ is the Killing form. The curvature tensor applied to elements $X,Y,Z,W\in {\lie{g}}$ is then 
$(R(X,Y)Z,W)=([X,Y],[Z,W])/4$
so the Ricci tensor $S$ is $-B(X,Y)/4$ and the scalar curvature
$s=\dim G/4$. Taking an orthonormal basis $\{y_a\}$ of ${\lie{g}}$ the above expression for the curvature is
$$R=\frac{1}{8} \sum \ad(y_a)\otimes \ad(y_a)\in \lie{so}(n)\otimes \lie{so}(n).$$
Thus from our point of view, the curvature term is 
$$-\frac{1}{2}\sum \rho(\ad(y_a))^2$$
and this involves the quadratic Casimir of the composition of the adjoint representation and the spin representation. We already know this is a scalar, $s/4$, so it is enough to take an irreducible component of the representation to calculate it. If $w$ is the highest weight of this representation and $\delta$ is half the sum of the positive roots of $G$ then the general formula for the Casimir is $\langle w, w+2\delta\rangle$. Now if $\pm x_i, (1\le i\le k=[n/2])$ are the weights of the vector representation of $SO(n)$, the highest weight of the spin representation is $(x_1+\dots +x_k)/2$. But restricted to $\lie{g}$, the $x_i$ are the roots. Hence $\delta$ is the highest weight of an irreducible component and we deduce that $\langle \delta, \delta+2\delta\rangle=3\langle \delta, \delta \rangle$ is equal to the scalar curvature term $-s/2=-\dim G/8$. In other words,
$$\dim G=24\Vert \delta\Vert^2$$
which is the strange formula of Freudenthal and de Vries. Stripped of the connection to curvature, this is the argument in \cite{FS}. 
\subsection{The Hodge Laplacian}
A natural generalization of the Dirac operator is to insert a coefficient bundle $E$ with connection. There is then a connection $\nabla$ on the tensor product $V\otimes E$ and a Dirac operator 
$$D=(C\otimes 1)(\nabla_V\otimes 1+1\otimes \nabla_E)$$
acting on sections of $V\otimes E$. Then 
$$D^2=\nabla^*\nabla+\frac{1}{4}s+C^2(1\otimes R_E)$$
where $R_E$, a section of $\End E\otimes \Lambda^2T^*$, is the curvature of $\nabla_E$. On a Riemannian manifold we can take $E$ to be associated to a representation $\sigma$ of $Spin(n)$ and its induced connection and then 
$$C^2(1\otimes R_E)=-4\sum R_{ab}\rho(x_a)\otimes \sigma(x_b)$$
and the curvature term in the Weitzenb\"ock decomposition is 
\begin{equation}
-4\sum R_{ab}(\rho(x_a) \rho(x_b)\otimes 1 +\rho(x_a)\otimes \sigma(x_b)).
\label{weitz}
\end{equation}
The operator $K$ for the tensor product of the two representations is 
$$\sum R_{ab}(\rho(x_a)\otimes 1 +1\otimes \sigma(x_a))(\rho(x_b)\otimes 1 +1\otimes \sigma(x_b)).$$
Now take $\sigma=\rho$, the spin representation,  then  this becomes, since $R_{ab}$ is symmetric 
$$\sum R_{ab}(\rho(x_a)\rho(x_b)\otimes 1 + 2\rho(x_a)\otimes \rho(x_b)+1\otimes \rho(x_a)\rho(x_b)) .$$
But $\sum R_{ab}\rho(x_a)\rho(x_b)$ is a scalar, a multiple of the scalar curvature $s$, so we obtain
$$K=2\sum R_{ab}(\rho(x_a) \rho(x_b)\otimes 1 +\rho(x_a)\otimes \rho(x_b))$$ and the spinor Laplacian is a Laplacian in the sense of Definition \ref{Lap} for $t=-2$.

On an even-dimensional manifold the bundle $V\otimes V\cong V\otimes V^*\cong \End V$ is the bundle of Clifford algebras which is isomorphic as a vector bundle with connection to the exterior algebra $\Lambda^*T^*$. Moreover, the Dirac operator $D$ with coefficients in $V$ is the  operator $d+d^*$ and $D^2= dd^*+d^*d$ is the Hodge Laplacian. So we have shown that the Hodge Laplacian is of our type. In odd dimensions $V\otimes V$ is half the exterior algebra but the same Laplacian results.

\begin{ex} The curvature term for one-forms consists of taking the defining vector representation of $SO(n)$ and in our formalism it is 
$$-2\sum R_{ab}x_ax_b$$
where $\{x_a\}$ is an orthonormal basis of skew-symmetric matrices. In terms of the Riemann tensor $R_{ijk\ell}$ this is 
$$S_{i\ell}=-2\frac{1}{2}\sum_j R_{ijj\ell}$$
which is the Ricci tensor and is the original use of this decomposition in \cite{Boc}. Taking $t=+2$ gives the Laplace equation satisfied by Killing vector fields.
\end{ex}

One of the immediate consequences of our formalism is that, via Hodge theory and the de Rham theorem, the cohomology of a  Riemannian manifold with special holonomy $H\subset SO(n)$ breaks up into pieces corresponding to the decomposition of $\Lambda^*T^*$ into irreducible representations of $H$. This leads, for example, to the commuting $SU(1,1)$ action on the cohomology of a K\"ahler manifold where $H=U(m)\subset SO(2m)$ and the $Sp(1,1)$-action for a hyperk\"ahler manifold where $H=Sp(m)\subset SO(4m)$.
\subsection{Vector bundle-valued forms}\label{vector}

A number of Laplacians arise in Riemannian geometry in the context of   differential forms with values in a vector bundle associated to the metric. For example, the Riemann curvature tensor itself is a section $R$ of $\Lambda^2T^*\otimes \Lambda^2T^*$ and satisfies the second Bianchi identity which can be written as $d_{\nabla}R=0$ where $d_{\nabla}:\Omega^p(\Lambda^2T^*)\rightarrow \Omega^{p+1}(\Lambda^2T^*)$ is the extension of the exterior derivative using the Levi-Civita connection $\nabla$. A metric is said to have {\it harmonic curvature} if $d^*_{\nabla}R=0$ using the formal adjoint $d^*_{\nabla}:\Omega^{p+1}(\Lambda^2T^*)\rightarrow \Omega^{p}(\Lambda^2T^*)$. The Levi-Civita connection of   such a metric satisfies the Yang-Mills equations. Examples include Einstein manifolds and conformally flat manifolds of constant scalar curvature. 

There is a Laplacian $d_{\nabla}d^*_{\nabla}+d^*_{\nabla}d_{\nabla}$ here and we shall show that it fits our pattern. The operator $D=d_{\nabla}+d^*_{\nabla}$ is a Dirac operator with coefficient bundle $V\otimes \Lambda^2T^*$ but since $d_{\nabla}^2$ is non-zero in this case (in fact is the curvature itself) to get the Laplacian $\Delta$ we must project orthogonally $D^2$ onto the subbundle $\Lambda^2T^*\otimes \Lambda^2T^*\subset V\otimes (V^{\otimes 3})$. Here is a general lemma which helps the identification:

\begin{lemma} Let $D$ be the Dirac operator on $V\otimes V^{\otimes {k-1}}=V^{\otimes {k}}$.  Let $E\subset V^{\otimes {k}}$ be a bundle defined by a representation  on  which a transitive  group  $\Gamma$ of  permutations of the $k$ factors acts trivially. Let $P_E$ denote orthogonal projection onto $E$. Then the operator $P_ED^2P_E$ is a Laplacian in the sense of Definition \ref{Lap} with $t=-4/k$. 
\end{lemma}

Before proving this, let us see how it applies to the above Laplacian, for simplicity in the even-dimensional case, where $E$ is the space of curvature tensors.  We have observed that $V$ is isomorphic to $V^*$ and this identification is either skew or symmetric depending on the dimension. In other words there is  an invariant bilinear form  on spinors preserved by $Spin(n)$.  The Lie algebra $\rho(\lie{so}(n))$ therefore lies in $\Sym^2 V$ if the form is skew and in $\Lambda^2 V$ if it is symmetric: in either case   transposition of the factors acts as $\pm 1$.
So $\Lambda^2T^*\otimes \Lambda^2T^*\subset (V\otimes V)\otimes (V\otimes V)$ is preserved by the permutation  $(12)(34)$. But the curvature tensor $R=\sum R_{ab}x_a\otimes x_b$ has $R_{ab}=R_{ba}$ so it is invariant by   $(13)(24)$. These elements generate a transitive group of permutations of $\{1,2,3,4\}$ and so from the lemma we have an admissible Laplacian with $t=-1$.

We now prove the lemma.

\begin{lemprf} The permutation action on $V^{\otimes {k}}$ commutes with the $Spin(n)$-action so $P_Eg=gP_E$ for $g\in \Gamma$. But the action on $E$ is trivial so this is $P_E$. From (\ref{weitz}) the curvature term in the Weitzenb\"ock decomposition for $D^2$ is 
$$W=-4\sum_{a,b} R_{ab}[\rho(x_a)\rho(x_b)\otimes 1^{k-1}  +\rho(x_a)\otimes (\rho(x_b)\otimes 1^{k-2}+1\otimes\rho(x_b)\otimes 1^{k-3}+\dots )]$$
while the natural curvature term in our definition is 
$$tK=t\sum_{a,b} R_{ab}(\rho(x_a)\otimes 1^{k-1} +1\otimes \rho(x_a)\otimes 1^{k-2}+ \dots)(\rho(x_b)\otimes 1^{k-1} +1\otimes \rho(x_b)\otimes 1^{k-2} +\dots).$$
Each of the $k$ terms in $W$ is transformed as $g^{-1}Wg$ by an element of the transitive group into one of the $k^2$ terms in $K$.   Since $P_Eg^{-1}WgP_E=P_EWP_E$ and $K$ respects the decomposition into invariant subspaces, it follows that 
$$K=P_EKP_E=-\frac{k}{4}P_EWP_E.$$
The covariant derivative commutes with $P_E$ and it follows that 
$$P_E D^2P_E=\nabla^*\nabla-\frac{4}{k}K.$$
\end{lemprf} 

Another example concerns the second fundamental form $h$ of a hypersurface, a section of $\Sym^2 T^*$. Regarding $h\in \Omega^1(T^*)$, it   satisfies the Codazzi equation $d_{\nabla}h=0$. For a minimal hypersurface in the sphere we also have $d^*_{\nabla}h=0$. Because of the symmetry of $h$ it satisfies the conditions of the lemma and so $h$ is annihilated by a Laplacian of our form. 

Berger and Ebin \cite{BE} consider a number of operators besides this one on symmetric tensors. One is the Lichnerowicz Laplacian whose nullspace  measures infinitesimal deformations of Ricci-flat metrics. This is one of a generalization to all tensors which Lichnerowicz made and corresponds to our curvature term with $t=-2$.
\section{Positivity of curvature}

The value of Weitzenb\"ock decompositions lies principally in vanishing theorems. In our situation if the curvature term $tK$ on a compact manifold is a {\it positive} self-adjoint transformation  then the null space of the Laplacian vanishes. If it is positive semi-definite then any solution to $\Delta \psi=0$ is covariant constant. This can be seen by   integration:
$$0=\int_M(\Delta \psi,\psi)=\int_M(\nabla^*\nabla\psi,\psi)+(tK\psi,\psi)=\int_M(\nabla\psi,\nabla\psi)+\int_M(tK\psi,\psi)\ge 0.$$
Thus for spinors, positive scalar curvature gives vanishing and for $1$-forms positive Ricci tensor, the original vanishing theorem of Bochner \cite{Boc}. 
Similarly, for a compact Lie group, we see immediately that the cohomology is represented by bi-invariant forms  since $K$ is the quadratic Casimir, which is positive semi-definite and zero only on the trivial representation. 

If the matrix $R_{ab}$ is positive definite at each point of a manifold then this is the concept of  {\it positive curvature operator}. We prove next:
\begin{prp} The operator $-K$ is positive for the Laplacians associated to all non-trivial irreducible representations of $SO(n)$ if and only if $M$ has   positive curvature operator.
\end{prp}
\begin{prf} If $R_{ab}$ is positive definite at each point then there is an orthonormal basis $\{x_a\}$ of ${\lie{so}(n)}$ such that, for any representation $\sigma$  
$$K=\sum_a\lambda_a\sigma(x_a)^2$$
with $\lambda_a>0$.
Each term $\sigma(x_a)$ is skew adjoint, so $\sigma(x_a)^2$ is negative semi-definite. If $Kv=0$  for some $v$ then $\sigma(x_a)v=0$ for all $a$ and hence the representation has a trivial factor and is reducible.  Hence $-K$ is positive.

Conversely, suppose $-K$ is positive for all $\sigma$ then at each point $p\in M$ the expression  $\sum R_{ab}(p)\sigma(x_a)\sigma(x_b)$ is negative. Considering each $x_a$ as a left-invariant vector field on the group manifold $SO(n)$ we have a second-order linear differential operator acting on functions on $SO(n)$
$$L=\sum R_{ab}(p)x_ax_b.$$
Since  $L^2(SO(n))$ decomposes as irreducible representations, $L$  
acting on functions  is essentially positive (positive on the orthogonal complement to the constant functions). It follows that   its principal symbol is negative. Hence $R_{ab}(p)$ is positive definite. 
\end{prf}

B\"ohm and Wilking showed by Ricci-flow techniques that positive curvature operator implies that $M$ is diffeomorphic to a spherical space form. It is immediate from our approach that if also the curvature is harmonic, then it is isometric to such a manifold, because we get  vanishing for all but the constant term in the decomposition of the curvature tensor into irreducible components. This is a  theorem of Tachibana \cite{Tac}.

\section{Postscript}
The reader may be curious to know what the other two open problems were and what became of them in the past 35 years. Here they are:
\begin{quotation}
{\bf Problem C-2.} Consider a connected sum $M=S^3\times S^1\#\cdots \#S^3\times S^1$. Schoen \& Yau have shown that $M$ admits a conformally flat metric of positive scalar curvature. This in itself makes $M$ a good subject to study from the point of view of self-duality since all the vanishing theorems which apply to instanton bundles on $S^4$ apply to $M$. However, it would be more interesting to know if $M$ has a conformally flat metric of {\it constant} positive scalar curvature.

The Riemannian connection of such a metric would satisfy the non-self-dual Yang-Mills equations and it is possible that by deforming the conformal structure one would obtain in the limit a singular solution which may relate to ``multi-meron" solutions of the Yang-Mills equations.
\end{quotation}
 
Given the discussion in the previous section, I should explain what ``all" meant in this context.  This is the vanishing of the null-space of the Dirac operator on anti-self-dual spinors with coefficient bundle self-dual, and vanishing for the Laplacian on anti-self-dual $2$-forms with the same coefficient bundle. The present formalism makes this obvious, the essential point being that in four dimensions the component of $R$ in $ \Lambda_+^2T^*\otimes \Lambda_-^2 T^*$ is the trace-free Ricci tensor. Since $ \Lambda_+^2T^*$ is the Lie algebra of the copy of $SU(2)\subset SU(2)\times SU(2)=Spin(4)$ which acts trivially on the anti-self-dual spinors and forms, the trace-free Ricci tensor does not contribute to the curvature term $K$. For a conformally flat manifold the only term left in the curvature is the scalar curvature, hence the vanishing theorem. 
As to the open question, Schoen's resolution of the Yamabe problem \cite{Sch} gave a positive solution.  

The second part relates to the observation that the classical meron solution of the Yang-Mills equations with a singularity at the origin in $\R^4$ \cite{Alf} can be interpreted as the standard metric of constant scalar curvature on $S^3\times S^1$,  as observed in the physics literature in \cite{Chu}. From the Riemannian point of view this is the observation in Section \ref{vector} that a conformally flat manifold with constant scalar curvature has a harmonic metric. In four dimensions the Yang-Mills equations are conformally invariant, and the  connected sum $M$ is the quotient of an open set in $\R^4$ by a discrete group of conformal transformations, so this metric   describes a Yang-Mills connection in an open set of Euclidean space. 

\begin{quotation}
{\bf Problem C-3.} What is the complex structure on the three-dimensional ``twistor space" $Z\rightarrow P_1$, obtained from the Calabi-Yau theorem for a K3 surface? From the Penrose point of view this is an easier problem, but a necessary prerequisite, for computing the metric explicitly.
\end{quotation}

In many ways this question is ill-formulated. In what terms  is one meant to describe a complex structure of a non-algebraic complex manifold? It seems as if numerical methods provide  the only approach to any sort of explicitness at the moment \cite{SKD}. However, vanishing theorems contribute to identifying one aspect of the problem: the only section of a holomorphic line bundle on $Z$ is the pull-back of  a section of ${\mathcal O}(k)$ on $P_1$ for some $k$. And without any sections it is difficult to find any holomorphic data to begin to describe $Z$ explicitly. 

To see this,  let $L$ be a holomorphic line bundle on $Z$. Restricted to each twistor line, which is a holomorphic section of $\pi:Z\rightarrow P_1$, it has degree $d$. Applying the real structure $\sigma$ on $Z$ we get a real holomorphic line bundle $L\otimes \sigma^*L$ and then $L\otimes\sigma^*L\otimes \pi^*{\mathcal O}(-2d)$ is real and trivial on each line. By the Atiyah-Ward correspondence it defines a real line bundle on the K3 surface $M$ with anti-self-dual connection. A real line bundle has structure group which reduces to $\{\pm 1\}$ and hence its real Chern class is zero.  But then we have a harmonic form, its curvature, which is  cohomologically trivial, a contradiction. It follows that $L\cong U\otimes \pi^*{\mathcal O}(d)$ where $\sigma^*U\cong U^*$ and is trivial on each line.
The crucial point now is that  $U$ is the pull-back of  an anti-self-dual  Hermitian line bundle on $M$ via the non-holomorphic projection from $Z=M\times S^2$ to $M$. 

Since $\pi^*{\mathcal O}(d)$ is holomorphically trivial on each fibre of $\pi$, if the line bundle $L$ had a holomorphic section on $Z$ then so would $U$ on $M$, over each point of $P_1$. But this vanishing theorem was precisely the one which motivated the Hitchin-Kobayashi conjecture.

 \end{document}